\documentclass[12pt,twoside]{article}

\usepackage[english]{babel}
\usepackage{umj_1}           

\usepackage{graphicx, srcltx}   
\usepackage{cite}               

\usepackage{amsmath}            
\usepackage{amsfonts,amssymb}   
\usepackage{srcltx}             

\newcommand{\Leftclt}{Roman A. Veprintsev}
\newcommand{\Rightclt}{Dunkl harmonic analysis and fundamental sets of functions on the sphere}
\usepackage{authblk}

\begin{document}

\begin{Titul}
{\large \bf DUNKL HARMONIC ANALYSIS\\ AND FUNDAMENTAL SETS OF FUNCTIONS\\[0.2em] ON THE UNIT SPHERE }\\[3ex]
{{\bf Roman~A.~Veprintsev} \\[5ex]}
\end{Titul}

\begin{Anot}
{\bf Abstract.} Using Dunkl theory, we introduce into consideration some weighted $L_p$-spaces on $[-1,1]$ and on the unit Euclidean sphere $\mathbb{S}^{d-1}$, $d\geq 2$. Then we define a family of linear bounded operators $\{V_\kappa^p(x)\colon x\in\mathbb{S}^{d-1}\}$ acting from the $L_p$-space on $[-1,1]$ to the $L_p$-space on $\mathbb{S}^{d-1}$, $1\leq p<\infty$. We establish a necessary and sufficient condition for a function $g$ belonging to the $L_p$-space on $[-1,1]$ such that the family of functions $\{V_\kappa^p(x;g)\colon x\in\mathbb{S}^{d-1}\}$ is fundamental in the $L_p$-space on $\mathbb{S}^{d-1}$.


\end{Anot}






\section{Introduction and preliminaries}

In this section we introduce some basic definitions and notions of general Dunkl theory (see, e.g., \cite{dai_xu_book_approximation_2013,dunkl_article_reflection_1988,dunkl_article_operators_1989,dunkl_article_integral_kernels_1991,dunkl_xu_book_orthogonal_polynomials_2014}); for a background on reflection groups and root systems the reader is referred to \cite{humphreys_book_reflection_groups_1990,dunkl_xu_book_orthogonal_polynomials_2014}.

Let $\mathbb{N}_0$ be the set of nonnegative integers, let $\mathbb{R}^d$ be the $d$-dimensional real Euclidean space of all $d$-tuples of real numbers. For $x\in\mathbb{R}^d$, we write $x=(x_1,\dots,x_d)$. The inner product of $x,\,y\in\mathbb{R}^d$ is denoted by $\langle x,y\rangle=\sum\limits_{i=1}^d x_iy_i$, and the norm of $x$ is denoted by $\|x\|=\sqrt{\langle x,x\rangle}$. Let $\mathbb{S}^{d-1}=\{x\colon\, \|x\|=1\}$ be the unit sphere in $\mathbb{R}^d$, $d\geq 2$. Denote by $d\omega$ the usual Lebesgue measure on $\mathbb{S}^{d-1}$.

For a nonzero vector $v\in\mathbb{R}^d$, define the reflection $s_v$ by
\begin{equation*}
s_v(x)=x-2\frac{\langle x,v\rangle}{\|v\|^2}\,v,\quad x\in\mathbb{R}^d.
\end{equation*}
Each reflection $s_v$ is contained in the orthogonal group $O(\mathbb{R}^d)$.

Recall that a finite set $R\subset\mathbb{R}^d\setminus\{0\}$ is called a root system if the following conditions are satisfied:

(1) $R\cap\mathbb{R}v=\{\pm v\}$\, for all $v\in R$;\qquad (2) $s_v(R)=R$\, for all $v\in R$.

\noindent The subgroup $G=G(R)\subset O(\mathbb{R}^d)$ which is generated by the reflections $\{s_v\colon\, v\in R\}$ is called the reflection group associated with $R$. It is known that the reflection group $G$ is finite and the set of reflections contained in $G$ is exactly $\{s_v\colon\, v\in R\}$.

Each root system $R$ can be written as a disjoint union $R=R_+\cup (-R_+)$, where $R_+$ and $-R_+$ are separated by a hyperplane through the origin. Such a set $R_+$ is called a positive subsystem. Its choice in not unique.

A nonnegative function $\kappa$ on a root system $R$ is called a multiplicity function if it is $G$-invariant, i.e. $\kappa(v)=\kappa(g(v))$ for all $v\in R$, $g\in G$.

Note that definitions given below do not depend on the special choice of $R_+$, thanks to the $G$-invariance of $\kappa$.

The Dunkl operators are defined by
\begin{equation*}
\mathcal{D}_i f(x)=\frac{\partial f(x)}{\partial x_i}+\sum\limits_{v\in R_+} \kappa(v)\frac{f(x)-f(s_v(x))}{\langle v,x\rangle}\, v_i,\quad 1\le i\le d.
\end{equation*}
In case $\kappa=0$, the Dunkl operators reduce to the corresponding partial derivatives. These operators were introduced and first studied by C.\,F.~Dunkl.

Let
\begin{equation*}\label{special_indeces}
\lambda_\kappa=\gamma_\kappa+\frac{d-2}{2},\qquad\gamma_\kappa=\sum\limits_{v\in R_+} \kappa(v),
\end{equation*}
let $w_\kappa$ denote the weight function on $\mathbb{S}^{d-1}$ defined by
\begin{equation*}\label{Dunkl_weight_function}
w_\kappa(x)=\prod\limits_{v\in R_+} |\langle v,x\rangle|^{2\kappa(v)},\quad x\in\mathbb{S}^{d-1}.
\end{equation*}
The weight function $w_\kappa$ is a positively homogeneous $G$-invariant function of degree $2\gamma_\kappa$. 
In case $\kappa=0$, $w_\kappa$ is identically equal to $1$.

Suppose $\Pi^d$ is the space of all polynomials in $d$ variables with complex coefficients, $\mathcal{P}_n^d$ is the subspace of homogeneous polynomials of degree $n\in\mathbb{N}_0$ in $d$ variables.

C.\,F.~Dunkl has proved in \cite{dunkl_article_integral_kernels_1991} that there exists a unique linear isomorphism $V_\kappa$ of $\Pi^d$ such that
\begin{equation*}
V_\kappa(\mathcal{P}_n^d)=\mathcal{P}_n^d,\,\,\, n\in\mathbb{N}_0,\qquad V_\kappa 1=1,\qquad \mathcal{D}_i V_\kappa=V_\kappa\frac{\partial}{\partial x_i},\,\,\, 1\leq i\leq d.
\end{equation*}
This operator is called the Dunkl intertwining operator. The operator $V_\kappa$ was studied by many mathematicians (for example, C.\,F.~Dunkl, M.~R\"{o}sler, K.~Trim\`{e}che, Y.~Xu). If $\kappa=0$, then $V_\kappa$ is the identity operator.

Throughout this paper, we assume that $p\in[1,\infty)$ and $\lambda_\kappa>0$. In particular, it follows that $\gamma_\kappa>0$ if $d=2$.

To explain our main result of the present paper, we need to introduce some weighted $L_p$-spaces and one family of linear operators.

Denote by $L_{\kappa,p}(\mathbb{S}^{d-1})$ the space of complex-valued Lebesgue measurable functions $f$ on $\mathbb{S}^{d-1}$ with finite norm
\begin{equation*}
\|f\|_{\kappa,p,\mathbb{S}^{d-1}}=\Bigl(\int\nolimits_{\mathbb{S}^{d-1}} |f(x)|^p \,d\sigma_\kappa(x)\Bigr)^{1/p},\qquad d\sigma_\kappa(x)=a_\kappa\,w_\kappa(x)\,d\omega(x),
\end{equation*}
where the normalizing constant $a_\kappa$ satisfies $a_\kappa\int\nolimits_{\mathbb{S}^{d-1}} w_\kappa\,d\omega=1$. The space $L_{\kappa,2}(\mathbb{S}^{d-1})$ is a complex Hilbert space with the inner product
\begin{equation*}\label{inner_product_for_space_on_sphere}
\langle f,h\rangle_{\kappa,\mathbb{S}^{d-1}}=\int\nolimits_{\mathbb{S}^{d-1}} f(x)\overline{h(x)} \,d\sigma_\kappa(x).
\end{equation*}

We also introduce the space $L_{\kappa,\infty}(\mathbb{S}^{d-1})$ composed of all complex-valued Lebesgue measurable functions defined on $\mathbb{S}^{d-1}$ which are $\sigma_\kappa$-measurable and $\sigma_\kappa$-essentially bounded. Because $w_\kappa$ is $\omega$-a.e. nonzero on $\mathbb{S}^{d-1}$, the above notions coincides with the one of $w_\kappa$-measurable and $w_\kappa$-essentially bounded function, respectively.

Let $\lambda>0$. Suppose $L_{p,\lambda}[-1,1]$ is the space of complex-valued Lebesgue measurable functions $g$ on the segment $[-1,1]$ with finite norm
\begin{equation*}
\|g\|_{\lambda,p,[-1,1]}=\Bigl(c_\lambda\int\nolimits_{-1}^1 |g(t)|^p\, (1-t^2)^{\lambda-1/2}\,dt\Bigr)^{1/p},\qquad c_\lambda=\Bigl(\int\nolimits_{-1}^1 (1-t^2)^{\lambda-1/2}\,dt\Bigr)^{-1}.
\end{equation*}

The Gegenbauer polynomials $C_n^\lambda$ (see, e.g., \cite[p.~302]{andrews_askey_roy_book_special_functions_1999}) are orthogonal with respect to the weight function $(1-t^2)^{\lambda-1/2}$. For a function $g\in L_{p,\lambda}[-1,1]$, its Gegenbauer expansion takes the form
\begin{equation}\label{Gegenbauer_expansion}
g(t)\sim \sum\limits_{n=0}^\infty \Lambda_{n,\lambda}(g) \frac{n+\lambda}{\lambda} C_n^\lambda(t)\quad\text{with}\quad \Lambda_{n,\lambda}(g)=\frac{c_\lambda}{C_n^\lambda(1)}\int\nolimits_{-1}^1 g(t) C_n^\lambda(t) (1-t^2)^{\lambda-1/2}\,dt,
\end{equation}
since $\|C_n^\lambda\|_{\lambda,2,[-1,1]}^2=C_n^\lambda(1)\lambda/(n+\lambda)$.

Theorem~13.17 in \cite{hewitt_stromberg_book_real_analysis_1965} says that $L_{\kappa,q}(\mathbb{S}^{d-1})\subset L_{\kappa,p}(\mathbb{S}^{d-1})$ and $L_{q,\lambda}[-1,1]\subset L_{p,\lambda}[-1,1]$ for $1\leq p<q<\infty$.

Let us now consider the operators
\begin{equation*}
V_\kappa^p(x)\colon L_{p,\lambda_\kappa}[-1,1]\to L_{\kappa,p}(\mathbb{S}^{d-1})
\end{equation*}
depending on $x\in\mathbb{S}^{d-1}$ which are defined by the following rule:
\begin{equation*}
V_\kappa^p(x;g,y)=V_\kappa\bigl[g(\langle x,\cdot\rangle)\bigr](y),\quad y\in\mathbb{S}^{d-1}\quad \bigl(g\in L_{p,\lambda_\kappa}[-1,1]\bigr).
\end{equation*}
They are linear bounded operators with operator norm $1$.
It follows from the properties of $V_\kappa$ known from the general Dunkl theory (see especially the proof of Theorem~7.4.2 in \cite{dai_xu_book_approximation_2013}) that these operators are well defined.

The main result of this paper is a generalization of Theorem~2.4 in \cite{menegatto_article_fundamental_sets_1998}. More precisely, we establish a necessary and sufficient condition for a function $g\in L_{p,\lambda_\kappa}[-1,1]$ such that the family of functions \begin{equation}\label{fundamental_set}\mathcal{M}_\kappa^p(g)=\{V_\kappa^p(x;g)\colon x\in\mathbb{S}^{d-1}\}\end{equation} is fundamental in the space $L_{\kappa,p}(\mathbb{S}^{d-1})$. This fact is formulated and proved in section~\ref{formulation_and_proof_of_the_main_result}.

Recall that a set $\mathcal{F}$ in a Banach space $\mathcal{E}$ is said to be fundamental if the linear span of $\mathcal{F}$ is dense in $\mathcal{E}$. To prove the main result, we use a consequence of the Hahn--Banach theorem related to fundamentality of sets in normed linear spaces. We include it as a separate lemma for convenience.

\begin{lemen}\label{consequence_of_the_Hahn--Banach_theorem}
Let $\mathcal{F}$ be a subset of a Banach space $\mathcal{E}$. In order that $\mathcal{F}$ be fundamental in $\mathcal{E}$, it is necessary and sufficient that $\mathcal{F}$ not be annihilated by a nonzero bounded linear functional on $\mathcal{E}$.
\end{lemen}

\section{Some facts of Dunkl harmonic analysis on the unit sphere}

The Dunkl Laplacian $\Delta_\kappa$ is defined by
\begin{equation*}
\Delta_\kappa=\sum\limits_{i=1}^d\mathcal{D}_i^2
\end{equation*}
and it plays the role similar to that of the ordinary Laplacian.
It reduces to the ordinary Laplacian provided that $\kappa=0$.

A $\kappa$-harmonic polynomial $P$ of degree $n\in\mathbb{N}_0$ in $d$ variables is a homogeneous polynomial $P\in\mathcal{P}_n^d$ such that $\Delta_\kappa P=0$. Its restriction to the unit sphere is called the $\kappa$-spherical harmonic of degree $n$ in $d$ variables. Denote by $\mathcal{A}_n^d(\kappa)$ the space of $\kappa$-spherical harmonics of degree $n$ in $d$ variables. The $\kappa$-spherical harmonics of different degrees turn out to be orthogonal with respect to the weighted inner product $\langle\cdot,\cdot\rangle_{\kappa,\mathbb{S}^{d-1}}$ \cite[Theorem~1.6]{dunkl_article_reflection_1988}.

Let $C(\mathbb{S}^{d-1})$ be the space of complex-valued continuous functions on $\mathbb{S}^{d-1}$.

\begin{lemen}\label{property_of_fundamentality}
The set $\bigcup\limits_{n=0}^\infty\mathcal{A}_n^d(\kappa)$ is fundamental in $C(\mathbb{S}^{d-1})$ and in $L_{\kappa,p}(\mathbb{S}^{d-1})$, $1\leq p<\infty$.
\end{lemen}

\proofen
Theorem~3.14 in \cite{rudin_book_analysis_1987} states that the space $C(\mathbb{S}^{d-1})$ is dense in $L_{\kappa,p}(\mathbb{S}^{d-1})$ for $1\leq p<\infty$. So it is sufficient to show that $\bigcup\limits_{n=0}^\infty\mathcal{A}_n^d(\kappa)$ is fundamental in $C(\mathbb{S}^{d-1})$.

By the Weierstrass approximation theorem, if $f$ is continuous on $\mathbb{S}^{d-1}$, then it can be uniformly approximated by polynomials restricted to $\mathbb{S}^{d-1}$. According to \cite[Theorem~1.7]{dunkl_article_reflection_1988}, these restrictions belong to the linear span of $\bigcup\limits_{n=0}^\infty\mathcal{A}_n^d(\kappa)$. Thus, $\bigcup\limits_{n=0}^\infty\mathcal{A}_n^d(\kappa)$ is fundamental in $C(\mathbb{S}^{d-1})$.
\hfill$\square$

The above proof is analogous to that of Corollary~2.3 in \cite{stein_weiss_book_Fourier_analysis_1971}.

\begin{lemen}
Let $g\in L_{p,\lambda_\kappa}[-1,1]$, $1\leq p<\infty$. Then for every $Y_n^\kappa\in\mathcal{A}_n^d(\kappa)$,
\begin{equation}\label{Funk-Hecke_formula}
\int\nolimits_{\mathbb{S}^{d-1}} V_\kappa^p(x;g,y)\,Y_n^\kappa(y) \,d\sigma_\kappa(y)=\Lambda_{n,\lambda_\kappa}(g)\,Y_n^\kappa(x),\quad x\in\mathbb{S}^{d-1},
\end{equation}
where the constant $\Lambda_{n,\lambda_\kappa}(g)$ is defined from \eqref{Gegenbauer_expansion}.
\end{lemen}

Equality~\eqref{Funk-Hecke_formula} is the Funk\,--\,Hecke formula for $\kappa$-spherical harmonics written in our setting and designations (cf. \cite[Theorem~7.2.7]{dai_xu_book_approximation_2013}, \cite[Theorem~2.1]{xu_article_Funk--Hecke_formula_2000}).

\section{Main result: proof and its consequence}\label{formulation_and_proof_of_the_main_result}

We can now state and prove the main theorem.

\begin{teoen}\label{main_result}
Let $d\geq 2$, $1\leq p<\infty$. Fix a root system $R$ in $\mathbb{R}^d$ and a multiplicity function $\kappa$ on $R$. Let $g\in L_{p,\lambda_\kappa}[-1,1]$. In order that the set $\mathcal{M}_\kappa^p(g)$ \eqref{fundamental_set} be fundamental in $L_{\kappa,p}(\mathbb{S}^{d-1})$, it is necessary and sufficient that $\Lambda_{n,\lambda_\kappa}(g)\not=0$ \eqref{Gegenbauer_expansion} for every $n\in\mathbb{N}_0$.
\end{teoen}

\proofen
We first prove that the condition is sufficient. Let $\Phi$ be a bounded linear functional on $L_{\kappa,p}(\mathbb{S}^{d-1})$ which annihilates $\mathcal{M}_\kappa^p(g)$. According to the Riesz representation theorem \cite[Theorem~6.16]{rudin_book_analysis_1987}, $\Phi$ can be written as follows: $\Phi(\,\cdot\,)=\langle \,\cdot\,,h\rangle_{\kappa,\mathbb{S}^{d-1}}$, where $h\hm\in L_{\kappa,q}(\mathbb{S}^{d-1})$ and $q$ is the exponent conjugate to $p$ ($p^{-1}+q^{-1}=1$; $q=\infty$ whenever $p=1$). Then the annihilating property of $\Phi$ reduces to
\begin{equation*}
\int\nolimits_{\mathbb{S}^{d-1}} V_\kappa^p(x;g,y)\, \overline{h(y)} \, d\sigma_\kappa(y)=0,\quad x\in\mathbb{S}^{d-1}.
\end{equation*}
Next, we multiply both sides of the previous equality by $Y_n^\kappa(x)\in\mathcal{A}_n^d(\kappa)$, $n\in\mathbb{N}_0$, and integrate the resulting expression with respect to the measure $d\sigma_\kappa$. H\"{o}lder's inequality implies that $V_\kappa^p(x;g,y) \, \overline{h(y)} \, Y_n^\kappa(x)$ is $\sigma_\kappa\times\sigma_\kappa$-integrable over $\mathbb{S}^{d-1}\times\mathbb{S}^{d-1}$, and hence, using the Fubini theorem to interchange the order of integration, we get
\begin{equation*}
\int\nolimits_{\mathbb{S}^{d-1}} \overline{h(y)}\, \biggl(\int\nolimits_{\mathbb{S}^{d-1}} V_\kappa^p(x;g,y)\,Y_n^\kappa(x) \,d\sigma_\kappa(x)\biggr) \,d\sigma_\kappa(y)=0.
\end{equation*}
Using the symmetric relation \cite[formula~(7)]{li_song_article_inversion_formulas_2009}
\begin{equation*}
V_\kappa\bigl[g(\langle x,\cdot\rangle)\bigr](y)=V_\kappa\bigl[g(\langle y,\cdot \rangle)\bigr](x) \quad\text{$\sigma_\kappa\times\sigma_\kappa$-a.e. on\, $\mathbb{S}^{d-1}\!\times\!\mathbb{S}^{d-1}$}
\end{equation*}
and the Funk\,--\,Hecke formula~\eqref{Funk-Hecke_formula}, we obtain
\begin{equation*}
\Lambda_{n,\lambda_\kappa}(g)\,\langle Y_n^\kappa,h\rangle_{\kappa,\mathbb{S}^{d-1}}=0,\qquad Y_n^\kappa\in\mathcal{A}_n^d(\kappa),\quad n\in\mathbb{N}_0.
\end{equation*}
It follows from the condition that
\begin{equation*}
\langle Y_n^\kappa,h\rangle_{\kappa,\mathbb{S}^{d-1}}=0,\qquad Y_n^\kappa\in\mathcal{A}_n^d(\kappa),\quad n\in\mathbb{N}_0.
\end{equation*}
Thus, $\Phi$ annihilates $\bigcup\limits_{n=0}^\infty \mathcal{A}_n^d(\kappa)$. By continuity of $\Phi$ and Lemma~\ref{property_of_fundamentality}, $\Phi=0$ on $L_{\kappa,p}(\mathbb{S}^{d-1})$. Therefore, the set $\mathcal{M}_\kappa^p(g)$ is fundamental in $L_{\kappa,p}(\mathbb{S}^{d-1})$ by Lemma~\ref{consequence_of_the_Hahn--Banach_theorem}.

Let us now prove that the condition described in the theorem is necessary. Assume, to reach a contradiction, that there exists an index $m\in\mathbb{N}_0$ such that $\Lambda_{m,\lambda_\kappa}(g)=0$. Select any nontrivial $\kappa$-spherical harmonic $Y_m^\kappa\in\mathcal{A}_m^d(\kappa)$ and consider a measure $\mu$ defined on the Lebesgue subsets $\mathcal{L}$ of $\mathbb{S}^{d-1}$ by the rule
\begin{equation*}
\mu(B)=\int\nolimits_{B} Y_m^\kappa(x)\,\sigma_\kappa(x),\quad B\in\mathcal{L}.
\end{equation*}
This measure is nontrivial by its definition.

Using the Funk\,--\,Hecke formula~\eqref{Funk-Hecke_formula}, we obtain
\begin{equation*}
\begin{split}
\int\nolimits_{\mathbb{S}^{d-1}} V_\kappa^p(x;g,y)\,d\mu(y)&=\int\nolimits_{\mathbb{S}^{d-1}} V_\kappa^p(x;g,y) \, Y_m^\kappa(y)\,d\sigma_\kappa(y)\\
&=\Lambda_{m,\lambda_\kappa}(g)\,Y_m^\kappa(x)=0,\quad x\in\mathbb{S}^{d-1}.
\end{split}
\end{equation*}
Thus, the nontrivial bounded linear functional $\Phi_1$ on $L_{\kappa,p}(\mathbb{S}^{d-1})$ given by $\Phi_1(f)=\int\nolimits_{\mathbb{S}^{d-1}} f\,d\mu$ annihilates $\mathcal{M}_\kappa^p(g)$. By Lemma~\ref{consequence_of_the_Hahn--Banach_theorem}, $\mathcal{M}_\kappa^p(g)$ is not fundamental in $L_{\kappa,p}(\mathbb{S}^{d-1})$. This contradicts our assumption.
\hfill$\square$

The above proof is exactly like that of~Theorem~2.4 in \cite{menegatto_article_fundamental_sets_1998}. Using the scheme of the proof of the theorem, one can prove the following result.

\begin{coren}
Let $d\geq 2$, $s\geq 1$, $1\leq p<\infty$. Fix a root system $R$ in $\mathbb{R}^d$ and a multiplicity function $\kappa$ on $R$. Let $g_1,\dots,g_s\in L_{p,\lambda_\kappa}[-1,1]$. In order that the set $\bigcup\limits_{i=1}^s\mathcal{M}_\kappa^p(g_i)$ be fundamental in $L_{\kappa,p}(\mathbb{S}^{d-1})$, it is necessary and sufficient that $\sum\limits_{i=1}^s|\Lambda_{n,\lambda_\kappa}(g_i)|\not=0$ for every $n\in\mathbb{N}_0$.
\end{coren}

\begin{Biblioen}
\bibitem{andrews_askey_roy_book_special_functions_1999}Andrews, G.E., Askey, R. and Roy, R., {\em Special functions}. Cambridge University Press, 1999.

\bibitem{dai_xu_book_approximation_2013}Dai, F. and Xu, Y., {\em Approximation theory and harmonic analysis on spheres and balls}. Springer, Berlin--New York, 2013.

\bibitem{dunkl_article_reflection_1988}Dunkl, C.F., Reflection groups and orthogonal polynomials on the sphere. {\em Math. Z.} 197: 33--60, 1988.

\bibitem{dunkl_article_operators_1989}Dunkl, C.F., Differential-difference operators associated to reflection groups. {\em Trans. Amer. Math. Soc.} 311(1): 167--183, 1989.

\bibitem{dunkl_article_integral_kernels_1991}Dunkl, C.F., Integral kernels with reflection group invariance. {\em Can. J. Math.} 43(6): 1213--1227, 1991.

\bibitem{dunkl_xu_book_orthogonal_polynomials_2014}Dunkl, C.F. and Xu, Y., {\em Orthogonal polynomials of several variables}. 2nd ed., Cambridge University Press, 2014.

\bibitem{hewitt_stromberg_book_real_analysis_1965}Hewitt, E. and Stromberg, K., {\em Real and abstract analysis}. Springer-Verlag, New York, 1965.

\bibitem{humphreys_book_reflection_groups_1990}Humphreys, J.E., {\em Reflection groups and Coxeter groups}. Cambridge University Press, 1990.

\bibitem{li_song_article_inversion_formulas_2009}Li, Zh., Song, F., Inversion formulas for the spherical Radon--Dunkl transform. {\em SIGMA.} 5: 025, 15 pages, 2009.

\bibitem{menegatto_article_fundamental_sets_1998}Menegatto, V.A., Fundamental sets of functions on spheres. {\em Methods Appl. Anal.} 5(4): 387--398, 1998.

\bibitem{rudin_book_analysis_1987}Rudin, W., {\em Real and complex analysis}. 3rd ed., McGraw-Hill, New York, 1987.

\bibitem{stein_weiss_book_Fourier_analysis_1971}Stein, E.M. and Weiss, G., {\em Introduction to Fourier analysis on Euclidean spaces}. Princeton University Press, Princeton, 1971.

\bibitem{xu_article_Funk--Hecke_formula_2000}Xu, Y., Funk\,--\,Hecke formula for orthogonal polynomials on spheres and on balls. {\em Bull. London Math. Soc.} 32: 447--457, 2000.
\end{Biblioen}

\noindent \textsc{Independent researcher, Uzlovaya, Russia}

\noindent \textit{E-mail address}: \textbf{veprintsevroma@gmail.com}

\end{document}